\documentclass{crm-article}

\usepackage{citehack}
\usepackage{algorithm}
\usepackage{algorithmic}

\def\ag#1{{\color{black}#1}}

\newcommand*{\norm}[1]{\left\lVert#1\right\rVert}
\newcommand{\argmin}{\mathop{\mathrm{argmin}}}
\newcommand{\argdualmin}{\mathop{\mathrm{argdual}}}
\newcommand{\leqarg}[1]{\ensuremath{\stackrel{\text{#1}}{\leq}}}
\newcommand{\eqarg}[1]{\ensuremath{\stackrel{\text{#1}}{=}}}
\newcommand*{\e}{\varepsilon}

\begin{document}

%\year=2018 % current year by default
\journalVol{10}

%выпуска
\journalNo{1}
\setcounter{page}{1}

% раздел журнала
\journalSection{Математические основы и численные методы моделирования}
\journalSectionEn{Mathematical modeling and numerical simulation}

% дата получения
\journalReceived{01.06.2016.}
%\journalReviewed{01.06.2016.}

%принято к публикации
\journalAccepted{01.06.2016.}

%!!!!!!!! \englishpaper раскомментировать в том случае, если текст статьи на английском языке
%\englishpaper

%!!!!!!!! \affiliationnoref раскомментировать, если автор один или все авторы из одной организации
%\affiliationnoref

%!!!!!!!! \emailnoref раскомментировать в том случае, если автор единственный
%\emailnoref

\UDC{519.85}
\title{Прямо-двойственный быстрый градиентный метод с моделью}
\titleeng{Primal-dual fast gradient method with a model}
\thanks{Работа была поддержана грантом РФФИ 18-31-20005 мол-а-вед в первой части и грантом РНФ 17-11-01027 во второй.}%если имеется
\thankseng{This work was supported by RFFI 18-31-20005 mol\_a\_ved in the first part of the work and by RSCF grant No. 17-11-01027 in the second part of the work.}

%автор - в формате \author{\firstname{И.\,И.}~\surname{Иванов}}
\author{\firstname{А.\,И.}~\surname{Тюрин}}
%автор - в формате \authorfull{Имя Отчество Фамилия}
\authorfull{Александр Игоревич Тюрин}
% автор на англ. - в формате \authoreng{\firstname{I.\,I.}~\surname{Ivanov}}
\authoreng{\firstname{A.\,I.}~\surname{Tyurin}}
%автор на англ. - в формате \authorfull{Firstname M. Surname}
\authorfulleng{Alexander I. Tyurin}
%вписать свою электронную почту
\email{alexandertiurin@gmail.com}
%организация - в формате \affiliation{Московский государственный университет,\protect\\ Россия, 141700, г. Москва, ул. Университетская, д. 9}
\affiliation{Национальный исследовательский университет «Высшая школа экономики»,\protect\\ Россия, 101000, г. Москва, ул. Мясницкая, д. 20}
%организация - в формате \affiliationeng{Moscow State Institute University, 9 University street, Moscow, 141700, Russia}
\affiliationeng{National Research University Higher School of Economics,\protect\\ 20 Myasnitskaya ulitsa, Moscow, 101000, Russia}

% повторите блок для каждого автора;
% если авторов несколько, и автоматическая расстановка сносок от фамилий 
% к организациям приводит к неправильным результатам, укажите правильный 
% вариант в квадраных скобках

% \author[1,2,3]{\firstname{А.\,В.}~\surname{Гасников}}
% \authorfull{Александр Владимирович Гасников}
% \authoreng{\firstname{A.\,V.}~\surname{Gasnikov}}
% \authorfulleng{Alexander V. Gasnikov}
% \email{gasnikov@yandex.ru}
% \affiliation[2]{Московский физико-технический институт,\protect\\ Россия, 141701, г. Долгопрудный, Институтский пер., д. 9}
% \affiliationeng{Moscow Institute of Physics and Technology,\protect\\ 9 Institute lane, Dolgoprudny, 141701, Russia}
% \affiliation[3]{Институт проблем передачи информации РАН,\protect\\ Россия, 127051, г. Москва, Б. Каретный пер., д. 9}
% \affiliationeng{Institute for Information Transmission Problems RAS,\protect\\ 9 B. Karetny lane, Moscow, 127051, Russia}

\begin{abstract}
В данной работе рассматривается возможность применения концепции $(\delta, L)$--модели функции для оптимизационных задач, в которых посредством решения прямой задачи имеется необходимость восстанавливать решение двойственной задачи. Концепция $(\delta, L)$--модели основана на концепции $(\delta,L)$--оракула, предложенная Деволдерор–Глинером–Нестеровов, при этом данные авторы предложили фукнционалы в оптимизационных задачах аппроксимировать сверху выпуклой параболой с некоторым аддитивным шумом $\delta$, таким образом, им удалось получить квадратичные верхние оценки с шумом даже для негладких функционалов. Концепция $(\delta, L)$--модели продолжает эту идею за счет того, что аппроксимация сверху делается не выпуклой параболой, а некоторым более сложным выпуклым функционалом. Возможность восстанавливать решение двойственной задачи хорошо себя зарекомендовала себя, так как во многих случаях в прямой задаче можно значительно быстрее находить решение, чем в двойственной. Отметим, что прямо--двойственные методы хорошо изучены, но при этом, как правило, каждый метод предлагается под конкретный класс задач. Нашей же целью является предложить метод, который бы включал в себя сразу различные методы. Это реализуется за счет использования концепции $(\delta, L)$--модели и адаптивной структуры наших методов. Таким образом, нам удалось получить прямой--двойственный адаптивный градиентный метод и быстрый градиентный метод с $(\delta, L)$--моделью и доказать оценки сходимости для них, причем для некоторых классов задач данные оценки являются оптимальными. Основная идея заключается в том, что нахождение двойственных решений происходит относительно оптимизационной задачи, которая аппроксимируют прямую с помощью концепции $(\delta, L)$--модели и имеет более простую структуру, поэтому находить двойственное решение у нее проще. Стоит отметить, что это происходит на каждом шаге работы оптимизационного метода, таким образом, реализуется принцип "разделяй и властвуй". 
\end{abstract}


\begin{abstracteng}
In this work we consider a possibility to use the conception of $(\delta, L)$--model of a function for optimization tasks, whereby solving a primal problem there is a necessity to recover a solution of a dual problem. The conception of $(\delta, L)$--model is based on the conception of $(\delta, L)$--oracle which was proposed by Devolder--Glineur--Nesterov, herewith the authors proposed approximate a function with an upper bound using a convex quadratic function with some additive noise $\delta$. They managed to get convex quadratic upper bounds with noise even for nonsmooth functions. The conception of $(\delta, L)$--model continues this idea by using instead of a convex quadratic function a more complex convex function in an upper bound. Possibility to recover the solution of a dual problem gives great benefits in different problems, for instance, in some cases, it is faster to find a solution in a primal problem than in a dual problem. Note that primal--dual methods are well studied, but usually each class of optimization problems has its own primal--dual method. Our goal is to develop a method which can find solutions in different classes of optimization problems. This is realized through the use of the conception of $(\delta, L)$--model and adaptive structure of our methods. Thereby, we developed primal--dual adaptive gradient method and fast gradient method with $(\delta, L)$--model and proved convergence rates of the methods, moreover, for some classes of optimization problems the rates are optimal. The main idea is the following: we find a dual solution to an approximation of a primal problem using the conception of $(\delta, L)$--model. It is much easier to find a solution to an approximated problem, however, we have to do it in each step of our method, thereby the principle of "divide and conquer" is realized.

\end{abstracteng}
\keywordeng{fast gradient method}
\keywordeng{model of the function}
\keywordeng{primal--dual method}

\maketitle

%Раздел
\paragraph{Введение}
Методы оптимизации играют большую роль в решении различных задач. Важным свойством некоторых оптимизационных методов является их прямо--двойственность \cite{anikin2017dual, boyd2004convex, nesterov2015complexity, nesterov2009primal},~--- это возможность восстанавливать достаточно эффективно решение двойственной задачи по прямой (или наоборот). Данный подход хорошо себя зарекомендовал в транспортных задачах \cite{baymurzina2019universal, gasnikov2018dualtransport, gasnikoveffectivnie}, задаче машинного обучения SVM и многих других \ag{\cite{gasnikov2017universal}}. В данной работе мы предлагаем прямой--двойственный адаптивный градиентный и быстрый градиентный метод, использующий концепцию $(\delta, L)$--модели функции \cite{gasnikov2017universal, tyurin2017fast}, которая, в свою очередь, основана на концепции $(\delta, L)$--оракула \cite{devolder2014first, devolder2013first, devolder2013intermediate, devolder2013exactness}. Как и в ранних работах по $(\delta, L)$--модели методы из текущей работы
%  могут применимы к большому множеству задач, в частности, можно показать, что данные 
включают в себя классический градиентный метод \cite{nesterov2010introductory}, универсальный метод \cite{nesterov2015universal}, метод Франк–Вульфа \cite{ben-tal2015lectures}, композитная оптизация \cite{nesterov2013gradient}. Более того, концепция $(\delta, L)$--модели позволяет решать эффективно достаточно нетривиальные постановки задач \cite{stonyakin2019gradient, tyurin2017fast, gasnikov2017universal}. Для многих из них предложенные нами методы являются оптимальными \cite{nemirovskiy1979slognost, tyurin2017fast}.

\paragraph{Прямо--двойственный метод}

Рассмотрим общую задачу оптимизации \cite{nesterov2010introductory, vasiliev2017methods}:
\begin{equation}
    \label{main_task}
    f(x) \rightarrow \min_{x \in Q}.
\end{equation}
Функция $f(x)$ определена на некотором множестве $Q$, которое принадлежит линейному пространству $\mathbb{R}^{n}$:
$$f(x):Q\rightarrow\mathbb{R}, \quad Q \subset \mathbb{R}^{n}.$$
Далее и везде будем считать, что функции $f(x)$ выпуклая и на множестве $Q$ имеет хотя бы одну точку минимума, принадлежащую множество $Q$. Более того, будем предполагать, что множество $Q$ имеет следующий вид:
	\begin{align*}
	Q = \{x~|~x \in \widetilde{Q},~f_i(x) \leq 0~ \forall i \in [1, m]\},
	\end{align*}
где для любого $i$ функция $f_i(x):\widetilde{Q}\rightarrow\mathbb{R}$ выпуклая функция, и множество $\widetilde{Q}$ является выпуклым. Введем следующее обозначение:
	\begin{align*}
	F(x) = [f_1(x), \dots, f_m(x)]^T,
	\end{align*}
таким образом, получаем следующую задачу оптимизации:
	\begin{align}
	\label{main_task_dop}
	f(x)  \rightarrow \min_{x \in \widetilde{Q},~F(x) \leq 0}.
	\end{align}
Далее нам понадобится понятие \textit{прокс--функции} и \textit{дивергенции Брэгмана} \cite{ben-tal2015lectures}:
\begin{fed}
Функция $d(x):Q \rightarrow \mathbb{R}$ называется прокс--функцией, если $d(x)$ непрерывно дифференцируемая на $\textnormal{int }Q$ и $d(x)$ является 1--сильно выпуклой относительно нормы $\norm{\cdot}$ на множестве $\textnormal{int }Q$.
\end{fed}
\begin{fed}
Функция
\begin{align}
V(x,y) = d(x) - d(y) - \langle\nabla d(y), x - y\rangle
\end{align}
называется дивергенцией Брэгмана, где $d(x)$~--- произвольная прокс--функция.
\end{fed}
Из 1--сильной выпуклости прокс--функции моментально следует \cite{ben-tal2015lectures}, что \begin{gather}\label{bregman_strong}
V(x,y) \geq \frac{1}{2}\norm{x - y}^2.
\end{gather}

Введем понятие $(\delta, L)$--модели функции:
\begin{fed}	
\label{model}
	Пусть функция $\psi_{\delta}(x, y)$ выпуклая на множестве $Q$ и выполняется условие $\psi_{\delta}(x, x) = 0$ для всех $x \in Q$. Будем говорить, что $\psi_{\delta}(x, y)$ есть $(\delta, L)$--модель функции $f$ в точке $y$ относительно нормы $\norm{\cdot}$, если
	для любого $x \in Q$ неравенство
	\begin{gather}
	\label{model_def}
	0 \le f(x) - (f_{\delta}(y) + \psi_{\delta}(x, y)) \le \frac{L}{2}\norm{x - y}^2 + \delta
	\end{gather}
	выполнено для некоторых $L, \delta > 0$.
\end{fed}
Данное определение было введено и ранее в работах \cite{gasnikov2017universal, tyurin2017fast, stonyakin2019gradient} и базируется на концепции $(\delta, L)$--оракула \cite{devolder2013first, devolder2014first, devolder2013exactness}.

Найдем двойственную задачу \cite{boyd2004convex} для задачи \eqref{main_task_dop}, для этого выпишем следующее равенство:
	\begin{align*}
	\min_{x \in \widetilde{Q},~F(x) \leq 0} f(x) = \min_{x \in \widetilde{Q}} \max_{z \in \mathbb{R}^{m}_+} [f(x) + \langle z, F(x)  \rangle],
	\end{align*}
где $\mathbb{R}^{m}_+ = \{x~|~x \in \mathbb{R}^{n},~x_i \geq 0~ \forall i \in [1, m]\}$. В силу слабой двойственности \cite{boyd2004convex} будет верно неравенство:
	\begin{align*}
	\min_{x \in \widetilde{Q},~F(x) \leq 0} f(x) \geq \max_{z \in \mathbb{R}^{m}_+} \min_{x \in \widetilde{Q}} [f(x) + \langle z, F(x)  \rangle],
	\end{align*}
Пусть
	\begin{align}
	\label{dual_function}
	g(z) = \max_{x \in \widetilde{Q}} [- f(x)-\langle z, F(x)  \rangle],
	\end{align}
тогда 
    \begin{align}
    \label{dual_ineq}
	\min_{x \in \widetilde{Q},~F(x) \leq 0} f(x) \geq - \min_{z \in \mathbb{R}^{m}_+} g(z).
	\end{align}
Выражение слева называется прямой задачей, а справа - двойственной задачей, определим ее отдельно:
	\begin{align}
	\label{dual_task_dop}
	g(z)  \rightarrow \min_{z \in \mathbb{R}^{m}_+}.
	\end{align}
Далее будем предполагать, что выполнены условия сильной двойственности \cite{boyd2004convex}, одним следствием этого является то, что неравенство в \eqref{dual_ineq} переходит в равенство. Для решений прямых и двойственных задач введем следующее обозначение.
\begin{fed}
Пусть $x_*$ произвольное решение прямой задачи
\begin{equation}
    \label{tmp_primal}
    p(x) \rightarrow \min_{x \in \widetilde{Q},~G(x) \leq 0}.
\end{equation}
Точка $z_*$ произвольное решение двойственной задачи
\begin{align*}
	h(z)  \rightarrow \min_{z \in \mathbb{R}^{m}_+},
\end{align*}
для \eqref{tmp_primal}, где $z$~---  это двойственные переменные соответствующие ограничениям $G(x) \leq 0$. Введем оператор \textit{$\argdualmin$}, зависящий от функции $p(x)$ и $G(x)$, и возвращающий $x_*$ и $z_*$:
	\begin{align*}
	(x_*, z_*) := \argdualmin_{x \in \widetilde{Q}}(p(x), G(x)).
	\end{align*}
\end{fed}
Пусть $x_*$ и $z_*$ произвольные решения прямой и двойственной задачи из \eqref{dual_ineq}, таким образом:
	\begin{align*}
	(x_*, z_*) := \argdualmin_{x \in \widetilde{Q}}(f(x), F(x)).
	\end{align*}
Прежде чем доказывать теорему рассмотрим алгоритм \ref{Alg1}. \ag{Этот алгоритм является комбинацией градиентого спуска в модельной общности из работы \cite{gasnikov2017universal} c прямо-двойственным субградиентным методом из работы \cite{nesterov2009primal}.} На вход алгори\ag{тм}у подается начальная точка $x_0$, произвольная константа $L_0 > 0$ и последовательность $\{\delta_k\}_{k\ge0}$. Будем предполагать далее, что для $\delta_k$ и точки $x_k$ всегда найдется некоторая константа $L_{k+1} > 0$ такая, что существует $(\delta_k, L_{k+1})$--модель в точке $x_k$. 
% Таким образом, мы всегда будем выходить из шага \ref{loop_state} алгоритма \ref{Alg1}. 
Будем также считать, что данное требование выполнено и для алгоритма \ref{Alg2}. Отметим еще, что $i_k$ в шаге \ref{loop_state} алгоритма \ref{Alg1} находится обычным перебором от $0$ до бесконечности, но из условия о существовании $(\delta_k, L_{k+1})$--модели в точке $x_k$ следует, что это этот процесс конечен, более того, несложно показать, что в среднем минимальное целое число $i_k$ для которого выполнено \eqref{exitLDL_G} равно $1$ \cite{nesterov2015universal}.

\begin{algorithm}
\caption{{Прямо--двойственный градиентный метод с моделью функции.}}
\label{Alg1}
\hspace*{\algorithmicindent} \textbf{Input: }$x_0$~--- начальная точка, $L_0 > 0$ и
$\{\delta_k\}_{k\ge0}$.
\begin{algorithmic}[1]
\STATE $A_0 := 0$
\FOR{$k \geq 0$}
\STATE \label{loop_state} Найти минимальное целое число $i_k\geq 0$ такое, что
\begin{equation}\label{exitLDL_G}
f_{\delta_k}(x_{k+1}) \leq f_{\delta_k}(x_{k}) + \psi_{\delta_k}(x_{k+1}, x_{k}) + \frac{L_{k+1}}{2}\norm{x_{k+1} - x_{k}}^2 + \delta_k,
\end{equation}
где $L_{k+1} := 2^{i_k-1}L_k$, $A_{k+1} := A_k + \frac{1}{L_{k+1}}$.
\begin{equation}\label{equmir2DL_G}
\phi_{k+1}(x) := \psi_{\delta_k}(x, x_k)+L_{k+1}V(x,x_k), \quad
(x_{k+1}, z_{k+1}) := \argdualmin_{x \in \widetilde{Q}}(\phi_{k+1}(x), F(x)).
\end{equation}
\ENDFOR
\end{algorithmic}
\hspace*{\algorithmicindent} \textbf{Output: } $\bar{x}_N= \frac{1}{A_N}\sum_{k=0}^{N-1}\frac{x_{k+1}}{L_{k+1}}$, $\bar{z}_N= \frac{1}{A_N}\sum_{k=0}^{N-1}\frac{z_{k+1}}{L_{k+1}}$
\end{algorithm}

\begin{lem}
	Пусть $\psi(x)$ выпуклая функция и 
	\begin{gather*}
	y = {\argmin_{x \in Q}}\{\psi(x) + V(x,u)\}.
	\end{gather*}
	Тогда выполнено неравенство
	\begin{equation*}
	\psi(x) + V(x,u) \geq \psi(y) + V(y,u) + V(x,y) \,\,\,\, \forall x \in Q.
	\end{equation*}
	\label{lemma_maxmin_2}
\end{lem}
Доказательство представлено в работе \cite{tyurin2017fast}, лемме 1.

\begin{cor}
	Пусть $\psi(x)$ выпуклая функция и 
	\begin{gather}
	\label{cor_1}
	(y, z) := \argdualmin_{x \in \widetilde{Q}}(\psi(x) + V(x,u), F(x)).
	\end{gather}
	Тогда выполнено неравенство
	\begin{equation*}
	\psi(x) + \langle z, F(x) \rangle + V(x,u) \geq \psi(y) + V(y,u) + V(x,y) \,\,\,\, \forall x \in \widetilde{Q}.
	\end{equation*}
	\label{cor_maxmin_2}
\end{cor}
\proof
Из \eqref{cor_1} и сильной двойственности следует, что 
	\begin{align}
	\label{cor_2}
	y = {\argmin_{x \in \widetilde{Q}}}\{\psi(x) + \langle z, F(x) \rangle + V(x,u)\}.
	\end{align}
Используя лемму \ref{lemma_maxmin_2} для \eqref{cor_2} получаем неравенство
    \begin{equation*}
	\psi(x) + \langle z, F(x) \rangle + V(x,u) \geq \psi(y) + \langle z, F(y) \rangle + V(y,u) + V(x,y) \,\,\,\, \forall x \in \widetilde{Q}.
	\end{equation*}
Из условия дополняющей нежёсткости \cite{boyd2004convex} верно, что $\langle z, F(y) \rangle = 0$. Следствие доказано.
\qed

Докажем теорему, которая дает оценки сходимости алгоритма \ref{Alg1}.

\begin{teo}
\label{teo_primal_dual}
	Пусть $x_0$~--- начальная точка, ($\bar{x}_N$, $\bar{z}_N$)~--- точки, полученные в результате работы алгоритма \ref{Alg1}, $x(\bar{z}_N)$~--- точка, в которой достигается максимум в \eqref{dual_function} при $z = \bar{z}_N$ и $V(x(\bar{z}_N)), x_0) \leq R^2,$ тогда будет верно неравенство
	\begin{align*}
	f(\bar{x}_N) \leq \min_{x \in \widetilde{Q}}\left[\frac{1}{A_N}\sum_{k=0}^{N-1}\frac{1}{L_{k+1}}(f_{\delta_k}(x_{k}) + \psi_{\delta_k}(x,x_{k})) + \langle \bar{z}_N, F(x) \rangle + \frac{V(x, x_0)}{A_N}\right] + \frac{1}{A_N}\sum_{k=0}^{N-1}\frac{2\delta_k}{L_{k+1}}
\end{align*}
и 
	\begin{align*}
f(\bar{x}_N) + g(\bar{z}_N) \leq \frac{R^2}{A_N} + \frac{1}{A_N}\sum_{k=0}^{N-1}\frac{2\delta_k}{L_{k+1}}.
\end{align*}
\end{teo}

\proof
Рассмотрим следующую цепочку неравенств:
	\begin{align*}
	f(x_{k+1}) &\leqarg{\eqref{model_def}}f_{\delta_k}(x_{k+1}) + \delta_k\\
	&\leqarg{\eqref{exitLDL_G}}f_{\delta_k}(x_{k}) + \psi_{\delta_k}(x_{k+1}, x_{k}) + \frac{L_{k+1}}{2}\norm{x_{k+1} - x_{k}}^2 + 2\delta_k\\
	&\leqarg{\eqref{bregman_strong}}f_{\delta_k}(x_{k}) + \psi_{\delta_k}(x_{k+1}, x_{k}) + L_{k+1}V(x_{k+1}, x_{k}) + 2\delta_k\\
	&\leqarg{С. \ref{cor_maxmin_2}}f_{\delta_k}(x_{k}) + \langle z_{k+1}, F(x) \rangle + \psi_{\delta_k}(x,x_{k})
	 	 + L_{k+1}V(x, x_k) - L_{k+1}V(x, x_{k+1}) + 2\delta_k.
	\end{align*}
Поделим неравенства на $L_{k+1}$, тогда получим
    \begin{align*}
	\frac{1}{L_{k+1}}f(x_{k+1}) \leq \frac{1}{L_{k+1}}(f_{\delta_k}(x_{k}) + \langle z_{k+1}, F(x) \rangle + \psi_{\delta_k}(x,x_{k}) + 2\delta_k)
	+ V(x, x_k) - V(x, x_{k+1}).
	\end{align*}
Если вычислить сумму неравенств по $k$ от $0$ до $N - 1$ и поделить на $A_N$, то будет верно
\begin{align*}
	\frac{1}{A_N}\sum_{k=0}^{N-1}\frac{1}{L_{k+1}}f(x_{k+1}) &\leq \frac{1}{A_N}\sum_{k=0}^{N-1}\frac{1}{L_{k+1}}(f_{\delta_k}(x_{k})+ \langle z_{k+1}, F(x) \rangle + \psi_{\delta_k}(x,x_{k}) + 2\delta_k)\\
	&+ \frac{1}{A_N}(V(x, x_0) - V(x, x_{N})).
	\end{align*}
Воспользуемся выпуклостью функции $f(x)$ и $V(x, x_{N}) \geq 0$, тогда получим
\begin{align*}
	f(\bar{x}_N) \leq \frac{1}{A_N}\sum_{k=0}^{N-1}\frac{1}{L_{k+1}}(f_{\delta_k}(x_{k}) + \langle z_{k+1}, F(x) \rangle + \psi_{\delta_k}(x,x_{k}) + 2\delta_k)
	+ \frac{V(x, x_0)}{A_N}.
\end{align*}
Данное неравенство верно для любого $x \in \widetilde{Q}$, поэтому
\begin{align*}
	f(\bar{x}_N) \leq \min_{x \in \widetilde{Q}}\left[\frac{1}{A_N}\sum_{k=0}^{N-1}\frac{1}{L_{k+1}}(f_{\delta_k}(x_{k}) + \psi_{\delta_k}(x,x_{k})) + \langle \bar{z}_N, F(x) \rangle + \frac{V(x, x_0)}{A_N}\right] + \frac{1}{A_N}\sum_{k=0}^{N-1}\frac{2\delta_k}{L_{k+1}}.
\end{align*}
Используя \eqref{model_def}, получаем неравенство 
\begin{align*}
	f(\bar{x}_N) \leq \min_{x \in \widetilde{Q}}\left[f(x) + \langle \bar{z}_N, F(x) \rangle + \frac{V(x, x_0)}{A_N}\right] + \frac{1}{A_N}\sum_{k=0}^{N-1}\frac{2\delta_k}{L_{k+1}}.
\end{align*}
По условию теоремы $x(\bar{z}_N) = \argmin_{x \in \widetilde{Q}} (-f(x) - \langle \bar{z}_N, F(x) \rangle)$, тогда по определению \eqref{dual_function} получаем цепочку неравенств
\begin{align*}
	f(\bar{x}_N) + g(\bar{z}_N) &\leq \frac{V(x(\bar{z}_N)), x_0)}{A_N} + \frac{1}{A_N}\sum_{k=0}^{N-1}\frac{2\delta_k}{L_{k+1}}\\
	&\leq \frac{R^2}{A_N} + \frac{1}{A_N}\sum_{k=0}^{N-1}\frac{2\delta_k}{L_{k+1}}.
\end{align*}
\qed

\paragraph{Быстрый прямо--двойственный метод}

В этом параграфе рассмотрим быстрый вариант градиентного метода. \ag{Этот алгоритм является комбинацией быстрого градиентого метода в модельной общности из работы \cite{tyurin2017fast} c прямо-двойственным субградиентным методом из работы \cite{nesterov2009primal}.} Сформулируем следующую теорему.

\begin{algorithm}
\caption{\bf{Быстрый прямо--двойственный градиентный метод с моделью функции}}
\label{Alg2}
\hspace*{\algorithmicindent} \textbf{Input: }$x_0$~--- начальная точка, $L_0 > 0$ и
$\{\delta_k\}_{k\ge0}$.
\begin{algorithmic}[1]
\STATE $y_0 := x_0$, $u_0 := x_0$, $\alpha_0 := 0$, $A_0 := \alpha_0$
\FOR{$k \geq 0$}
\STATE Найти минимальное целое число $i_k\geq 0$ такое, что
\begin{equation}
\begin{gathered}
f_{\delta_k}(x_{k+1}) \leq f_{\delta_k}(y_{k+1}) + \psi_{\delta_k}(x_{k+1}, y_{k+1}) +\frac{L_{k+1}}{2}\norm{x_{k+1} - y_{k+1}}^2 + \delta_k,
\label{exitLDL}
\end{gathered}
\end{equation}
где $L_{k+1}=2^{i_k-1}L_k$, 
\begin{equation}
\label{alpha_def}
\alpha_{k+1} := \frac{1 + \sqrt{ 1 + 4L_{k+1}A_{k}}}{2L_{k+1}}, \quad A_{k+1} := A_k + \alpha_{k+1}.
\end{equation}
\begin{gather}
y_{k+1} := \frac{\alpha_{k+1}u_k + A_k x_k}{A_{k+1}} \label{eqymir2DL}
\end{gather}
\begin{equation}\label{equmir2DL}
\phi_{k+1}(x)=V(x, u_k) + \alpha_{k+1}\psi_{\delta_k}(x, y_{k+1}),\quad
(u_{k+1}, z_{k+1}) := \argdualmin_{x \in \widetilde{Q}}(\phi_{k+1}(x), F(x))
\end{equation}
\begin{gather}
x_{k+1} := \frac{\alpha_{k+1}u_{k+1} + A_k x_k}{A_{k+1}} \label{eqxmir2DL}
\end{gather}
% \STATE \textbf{If} condition

% holds, \textbf{then}
% \begin{gather}
% L_{k+2} := \frac{L_{k+1}}{2}
% \end{gather}
% and go to the next step,
% \STATE \textbf{otherwise}
% \begin{gather}
% L_{k+1} := 2L_{k+1}
% \end{gather}
% \STATE and \textbf{repeat} the current step.
\ENDFOR
\end{algorithmic}
\hspace*{\algorithmicindent} \textbf{Output: } $x_N$, $\bar{z}_N= \frac{1}{A_N}\sum_{k=0}^{N-1}\alpha_{k+1}z_{k+1}$
\end{algorithm}

\begin{teo}
\label{fast_teo_primal_dual}
	Пусть $x_0$~--- начальная точка, ($x_N$, $\bar{z}_N$)~--- точки, полученные в результате работы алгоритма \ref{Alg2}, $x(\bar{z}_N)$~--- точка, в которой достигается максимум в \eqref{dual_function} при $z = \bar{z}_N$ и $V(x(\bar{z}_N)), x_0) \leq R^2,$ тогда будет верно неравенство
	\begin{align*}
	f(x_{N})\leq \min_{x \in \widetilde{Q}}\left[\frac{1}{A_N}\sum_{k=0}^{N-1}\alpha_{k+1}\Big(f_{\delta_k}(y_{k+1}) + \psi_{\delta_k}(x,y_{k+1})\Big) + \langle \bar{z}_N, F(x) \rangle + \frac{V(x, u_0)}{A_N}\right] + \frac{2}{A_N}\sum_{k=0}^{N-1}A_{k+1}\delta_k
\end{align*}
и 
	\begin{align*}
f(x_{N}) + g(\bar{z}_N) &\leq \frac{R^2}{A_N} + \frac{2}{A_N}\sum_{k=0}^{N-1}A_{k+1}\delta_k.
\end{align*}
\end{teo}

\proof
Рассмотрим следующую цепочку неравенств:
	\begin{align*}
	f(x_{k+1}) &\leqarg{\eqref{model_def}}f_{\delta_k}(x_{k+1}) + \delta_k\\
	&\leqarg{\eqref{exitLDL}} f_{\delta_k}(y_{k+1}) + \psi_{\delta_k}(x_{k+1},y_{k+1}) + \frac{L_{k+1}}{2}\norm{x_{k+1} - y_{k+1}}^2 + 2\delta_k\\
	&\eqarg{\eqref{eqxmir2DL}}
	f_{\delta_k}(y_{k+1}) + \psi_{\delta_k}\left(\frac{\alpha_{k+1}u_{k+1} + A_k x_k}{A_{k+1}},y_{k+1}\right) + \frac{L_{k+1}}{2}\norm{\frac{\alpha_{k+1}u_{k+1} + A_k x_k}{A_{k+1}} - y_{k+1}}^2 + 2\delta_k\\
    &\hspace{-20pt}\leqarg{\eqref{model_def}, \eqref{eqymir2DL}} f_{\delta_k}(y_{k+1}) +
	\frac{\alpha_{k+1}}{A_{k+1}}\psi_{\delta_k}(u_{k+1}, y_{k+1}) +\frac{A_k}{A_{k+1}}\psi_{\delta_k}(x_k, y_{k+1}) + \frac{L_{k+1} \alpha^2_{k+1}}{2 A^2_{k+1}}\norm{u_{k+1} - u_k}^2 + 2\delta_k\\
	 &=\frac{A_k}{A_{k+1}}(f_{\delta_k}(y_{k+1}) + \psi_{\delta_k}(x_k, y_{k+1}))+\frac{\alpha_{k+1}}{A_{k+1}}(f_{\delta_k}(y_{k+1}) +
	 \psi_{\delta_k}(u_{k+1}, y_{k+1}))\\
	 &+ \frac{L_{k+1} \alpha^2_{k+1}}{2 A^2_{k+1}}\norm{u_{k+1} - u_k}^2 + 2\delta_k
	 \end{align*}
	 Из \eqref{alpha_def} следует, что $A_{k+1} = L_{k+1}\alpha_{k+1}^2$, поэтому
	 \begin{align*}
	 f(x_{k+1})&=\frac{A_k}{A_{k+1}}(f_{\delta_k}(y_{k+1}) + \psi_{\delta_k}(x_k,y_{k+1}))\\
	 &+\frac{\alpha_{k+1}}{A_{k+1}}(f_{\delta_k}(y_{k+1}) + \psi_{\delta_k}(u_{k+1},y_{k+1})
	 + \frac{1}{2 \alpha_{k+1}}\norm{u_{k+1} - u_k}^2) + 2\delta_k\\
	 &\leqarg{\eqref{bregman_strong}} \frac{A_k}{A_{k+1}}(f_{\delta_k}(y_{k+1}) + \psi_{\delta_k}(x_k,y_{k+1}))\\
	 &+\frac{\alpha_{k+1}}{A_{k+1}}(f_{\delta_k}(y_{k+1}) + \psi_{\delta_k}(u_{k+1},y_{k+1})
	 + \frac{1}{\alpha_{k+1}}V(u_{k+1}, u_k)) + 2\delta_k\\
	 &\hspace{-10pt}\leqarg{С. \ref{cor_maxmin_2}, \eqref{model_def}} \frac{A_k}{A_{k+1}} f(x_k) +	 \frac{\alpha_{k+1}}{A_{k+1}}\Big(f_{\delta_k}(y_{k+1})+ \langle z_{k+1}, F(x) \rangle + \psi_{\delta_k}(x,y_{k+1}) \\
	 &+ \frac{1}{\alpha_{k+1}}V(x, u_k) - \frac{1}{\alpha_{k+1}}V(x, u_{k+1})\Big) + 2\delta_k.
	\end{align*}
	Умножим неравенство
	 \begin{align*}
	 f(x_{k+1})&\leq \frac{A_k}{A_{k+1}} f(x_k)	+ \frac{\alpha_{k+1}}{A_{k+1}}\Big(f_{\delta_k}(y_{k+1}) + \langle z_{k+1}, F(x) \rangle + \psi_{\delta_k}(x,y_{k+1})\\
	 &+ \frac{1}{\alpha_{k+1}}V(x, u_k) - \frac{1}{\alpha_{k+1}}V(x, u_{k+1})\Big) + 2\delta_k
	\end{align*}
	на $A_{k+1}$, тогда верно
	\begin{align*}
	 A_{k+1}f(x_{k+1})&\leq A_k f(x_k) +	\alpha_{k+1}\Big(f_{\delta_k}(y_{k+1}) + \langle z_{k+1}, F(x) \rangle + \psi_{\delta_k}(x,y_{k+1})\Big)\\
	 &+ V(x, u_k) - V(x, u_{k+1}) + 2A_{k+1}\delta_k.
	\end{align*}
	Если вычислить сумму неравенств по $k$ от $0$ до $N - 1$ и поделить на $A_N$, то будет верно
	\begin{align*}
	 f(x_{N})&\leq \frac{1}{A_N}\sum_{k=0}^{N-1}\alpha_{k+1}\Big(f_{\delta_k}(y_{k+1}) + \langle z_{k+1}, F(x) \rangle + \psi_{\delta_k}(x,y_{k+1})\Big)\\
	 &+\frac{1}{A_N}(V(x, u_0) - V(x, u_N)) + \frac{2}{A_N}\sum_{k=0}^{N-1}A_{k+1}\delta_k.
	\end{align*}
	Данное неравенство верно для любого $x \in \widetilde{Q}$. Кроме того, $V(x, u_{N}) \geq 0$, поэтому
	\begin{align*}
	 f(x_{N})\leq \min_{x \in \widetilde{Q}}\left[\frac{1}{A_N}\sum_{k=0}^{N-1}\alpha_{k+1}\Big(f_{\delta_k}(y_{k+1}) + \psi_{\delta_k}(x,y_{k+1})\Big) + \langle \bar{z}_N, F(x) \rangle + \frac{V(x, u_0)}{A_N}\right] + \frac{2}{A_N}\sum_{k=0}^{N-1}A_{k+1}\delta_k.
	\end{align*}
	Из \eqref{model_def} получаем неравенство 
	\begin{align*}
	 f(x_{N})\leq \min_{x \in \widetilde{Q}}\left[f(x) + \langle \bar{z}_N, F(x) \rangle + \frac{V(x, u_0)}{A_N}\right] + \frac{2}{A_N}\sum_{k=0}^{N-1}A_{k+1}\delta_k.
	\end{align*}
	По условию теоремы $x(\bar{z}_N) = \argmin_{x \in \widetilde{Q}} (-f(x) - \langle \bar{z}_N, F(x) \rangle)$, тогда по определению \eqref{dual_function} получаем цепочку неравенств
	\begin{align*}
	 f(x_{N}) + g(\bar{z}_N) &\leq \frac{V(x(\bar{z}_N), u_0)}{A_N} + \frac{2}{A_N}\sum_{k=0}^{N-1}A_{k+1}\delta_k\\
	 &\leq \frac{R^2}{A_N} + \frac{2}{A_N}\sum_{k=0}^{N-1}A_{k+1}\delta_k.
	\end{align*}
\qed

\paragraph{Анализ полученных результатов}
% Рассмотрим следующую задачу оптимизации:
% 	\begin{align}
% 	\label{dual_task}
% 	g(y)  \rightarrow \min_{y \in \widetilde{Q},~Ay \leq b}, 
% 	\end{align}
% где $g(y):\widetilde{Q}\rightarrow\mathbb{R}$, $\widetilde{Q} \in \mathbb{R}^{m}$, $A \in \mathbb{R}^{n, m}$ и $b \in \mathbb{R}^{n}$. Будем считать, что для $g(y)$ и $\widetilde{Q}$ выполнены все те же условия, что и для  $f(x)$ и $Q$ из задачи \eqref{main_task}.
% Целью данного параграфа является показать, как можно, используя алгоритм \ref{Alg1} и  \ref{Alg2}, найти решение не только прямой, но двойственной задачи \cite{boyd2004convex}.
% Рассмотрим задачу \eqref{main_task}.
% Дополнительно будет предполагать, что для задачи \eqref{dual_task} выполнены условия для сильной двойственности \cite{boyd2004convex}.

%   Ограничение подобного вида позволит нам далее построить следующую двойственную задачу \cite{boyd2004convex} для задачи \eqref{main_task}:
% 	\begin{align*}
% 	\label{dual_task}
% 	g(y) = \max_{x \in \widetilde{Q}} [\langle y , b - Ax\rangle - f(x)] \rightarrow \min_{y \in \mathbb{R}^{m}_+}, 
% 	\end{align*}
% где $\mathbb{R}^{m}_+ = \{y~|~y \in \mathbb{R}^{m},~y_i \geq 0~ \forall i \in [1, m]\}$.
% Выпишем неравенство из теоремы \ref{fast_teo_primal_dual} и запишем его с учетом условий на $Q$, получим
% 	\begin{align*}
% 	 f(x_{N})\leq \frac{1}{A_N}\min_{x \in \widetilde{Q},~Ax \leq b}\left[\sum_{k=0}^{N-1}\alpha_{k+1}\Big(f_{\delta_k}(y_{k+1}) + \psi_{\delta_k}(x,y_{k+1})\Big) + V(x, u_0)\right] + \frac{2}{A_N}\sum_{k=0}^{N-1}A_{k+1}\delta_k.
% 	\end{align*}
Для простоты будем далее считать, что $\delta_k \leq \e / 2$ в случае алгоритма \ref{Alg1} и $\delta_k \leq \e \alpha_{k+1} / (2A_{k+1})$ в случае алгоритма  \ref{Alg2}. Таким образом, мы получим, что слагаемое $$\frac{1}{A_N}\sum_{k=0}^{N-1}\frac{2\delta_k}{L_{k+1}}$$ из теоремы \ref{teo_primal_dual} и слагаемое $$\frac{2}{A_N}\sum_{k=0}^{N-1}A_{k+1}\delta_k$$ из теоремы \ref{fast_teo_primal_dual} будут меньше или равны $\e$.
Стоит отметить, что подобное ограничение на $\delta_k$ несильно сужает класс задач, в частности, в универсальных методах \cite{nesterov2015universal, gasnikov2017universal, tyurin2017fast} $\delta_k$ имеет именно такие порядки.

В работе \cite{tyurin2017fast} были получены различные оценки на $A_N$ в зависимости от дополнительной информации о функции $f(x)$, в данной работе все они сохраняются. В частности, когда $f(x)$ гладкая с $L$--липшицевым градиентом в норме $\norm{\cdot}$, в случае теоремы \ref{teo_primal_dual} $$A_N \geq \frac{2L}{N},$$
а в случае теоремы \ref{fast_teo_primal_dual}
$$A_N \geq \frac{8L}{(N+1)^2}.$$ Таким образом, алгоритмы \ref{Alg1} и \ref{Alg2} представляют из себя градиентный и быстрый градиентный метод соответственно. Отметим, что благодаря адаптивной структуре полученных методов они также являются универсальными \cite{nesterov2015universal, gasnikov2017universal, tyurin2017fast}. Более подробная информация о порядках $A_N$ содержится в работе \cite{tyurin2017fast}.

В теореме \ref{fast_teo_primal_dual} была получена скорости сходимости зазора двойственности $f(x_{N}) + g(\bar{z}_N)$. Если взять $N$ таким, чтобы $R^2/A_N \leq \e$, то мы получим неравенство \ag{$$0\le f(x_{N}) + g(\bar{z}_N) \leq 2\e.$$} Из неравенства \eqref{dual_ineq} мы получаем, что верны два других неравенства:
\ag{$$0\le f(x_{N}) - f(x_*) \leq 2\e, \quad 0\le g(\bar{z}_N) - g(z_*) \leq 2\e,$$} где $x_*$~--- оптимальное решение \eqref{main_task_dop}, а $z_*$~--- оптимальное решение \eqref{dual_task_dop}. Таким образом, нам удалось доказать, что алгоритм \ref{Alg2} генерирует $\e$--решение как для прямой, так и для двойственной задачи.
Все рассуждения о зазоре двойственности верны и для теоремы \ref{teo_primal_dual}.

% Кроме того, введем двойственные переменные $y \in \mathbb{R}^{m}_+$. Тогда
% 	\begin{align*}
% 	 f(x_{N})&\leq \frac{1}{A_N}\min_{x \in \widetilde{Q},~Ax \leq b}\left[\sum_{k=0}^{N-1}\alpha_{k+1}\Big(f_{\delta_k}(y_{k+1}) + \psi_{\delta_k}(x,y_{k+1})\Big) + V(x, u_0)\right] + \e \\
% 	 &\leq \min_{x \in \widetilde{Q}}\left[\frac{1}{A_N}\sum_{k=0}^{N-1}\alpha_{k+1}\Big(f_{\delta_k}(y_{k+1}) + \psi_{\delta_k}(x,y_{k+1})\Big) + \langle y \rangle + \frac{V(x, u_0)}{A_N}\right] + \e
% 	\end{align*}

%%%%%%%%%%%%%%%%%%
% Команды для оформления текста статьи, включая теоремы, следствия и пр. смотри в файле правила_оформления.txt

%% Список литературы
%пример оформления см. в файле правила_офорлмения.txt

\paragraph{Заключение}
В данный работе представлены два метода: прямо--двойственный градиентный и быстрый градиентный метод, и доказаны их оценки скорости сходимости. Данные методы являются обобщением алгоритмов, представленных в работе \cite{tyurin2017fast}. Стоит отметить, что все примеры задач из \cite{tyurin2017fast,stonyakin2019gradient} также применимы и здесь. Более того, данные алгоритмы применимы для задач, где важно одновременно находить решение прямой и двойственной задачи, решая только одну из этих задач, например, это используется в транспортных задачах \cite{baymurzina2019universal, gasnikov2018dualtransport, gasnikoveffectivnie}.

\end{document}